\documentclass[12]{article}
\usepackage{graphicx}
\usepackage[top=1in, bottom=1in, left=1in, right=1in]{geometry}

\usepackage{color}
\usepackage{amsmath}
\usepackage{relsize}

\begin{document}
\title{On the Location of the Non-Trivial Zeros of the\\ Riemann Hypothesis via Extended Analytic Continuation}
\author{Michael P. May}
\date{August 20, 2016}
\maketitle

\begin{abstract}
\begin{center}
This paper presents the results of a study on the analytic continuation of the prime zeta function to locate the non-trivial zeros of the Riemann Hypothesis in the critical strip of the complex plane.  A validation of RH is also proposed via the convergence and divergence of the $\zeta$, $m$, and $b$ functions.\\
\end{center}
\end{abstract}

\noindent One will recall the Euler-Riemann zeta function \cite{First50MillionPrimes}
\cite{HowEulerDiscoveredthezetafunction}:

\begin{equation}
\zeta(s) = 1 +\frac{1}{2^s}+\frac{1}{3^s}+\frac{1}{4^s}+\frac{1}{5^s}+ \cdots = \frac{1}{1-\frac{1}{2^s}}\cdot\frac{1}{1-\frac{1}{3^s}}\cdot\frac{1}{1-\frac{1}{5^s}}\cdot\frac{1}{1-\frac{1}{7^s}}\cdot\cdots
\end{equation}

\noindent I now introduce an equation that relates the prime sum and prime product functions \cite{A Non-Sieving Application of the Euler Zeta Function}:

\newcommand{\opM}{\mathop{\vphantom{\sum}\mathchoice
  {\vcenter{\hbox{\huge M}}}
  {\vcenter{\hbox{\Large A}}}{\mathrm{A}}{\mathrm{A}}}\displaylimits}

\begin{equation}
\displaystyle 1 + \displaystyle\sum\limits_{k=1}^\infty \frac{1}{p_k^s} = \frac{1+\displaystyle\sum\limits_{i=1}^\infty (-1)^{i}\displaystyle\sum\limits_{p}^{}{}_{i+1}}{\displaystyle\prod\limits_{k=1}^\infty {1-\frac{1}{\displaystyle p_k^s}•}}.
\end{equation}

\noindent The remainder function 

\begin{equation}
\displaystyle\sum\limits_{i=1}^\infty (-1)^{i}\displaystyle\sum\limits_{p}^{}{}_{i+1}
\end{equation}

\noindent found in the numerator of the right-hand side of Eq. 2 is defined here as the infinite alternating series

$$\--\Sigma_1 \Sigma_2+\Sigma_1 \Sigma_2 \Sigma_3-\Sigma_1 \Sigma_2 \Sigma_3 \Sigma_4+\Sigma_1 \Sigma_2 \Sigma_3 \Sigma_4 \Sigma_5-\Sigma_1 \Sigma_2 \Sigma_3 \Sigma_4 \Sigma_5 \Sigma_6+... $$

\noindent where

$$\Sigma_1 \Sigma_2 = \sum\limits_{i=1}^\infty \sum\limits_{j=i}^\infty \frac{1}{p_i^s} \frac{1}{p_j^s}$$

$$\Sigma_1 \Sigma_2 \Sigma_3 = \sum\limits_{i=1}^\infty \sum\limits_{j=i+1}^\infty \sum\limits_{k=i}^\infty \frac{1}{p_i^s} \frac{1}{p_j^s} \frac{1}{p_k^s}$$

$$\Sigma_1 \Sigma_2 \Sigma_3 \Sigma_4 = \sum\limits_{i=1}^\infty \sum\limits_{j=i+1}^\infty \sum\limits_{k=j+1}^\infty \sum\limits_{l=i}^\infty \frac{1}{p_i^s} \frac{1}{p_j^s} \frac{1}{p_k^s} \frac{1}{p_l^s}$$

$$\Sigma_1 \Sigma_2 \Sigma_3 \Sigma_4 \Sigma_5 = \sum\limits_{i=1}^\infty \sum\limits_{j=i+1}^\infty \sum\limits_{k=j+1}^\infty \sum\limits_{l=k+1}^\infty \sum\limits_{m=i}^\infty \frac{1}{p_i^s} \frac{1}{p_j^s} \frac{1}{p_k^s} \frac{1}{p_l^s}\frac{1}{p_m^s}$$

$$\Sigma_1 \Sigma_2 \Sigma_3 \Sigma_4 \Sigma_5 \Sigma_6= \sum\limits_{i=1}^\infty \sum\limits_{j=i+1}^\infty \sum\limits_{k=j+1}^\infty \sum\limits_{l=k+1}^\infty \sum\limits_{m=l+1}^\infty \sum\limits_{n=i}^\infty \frac{1}{p_i^s} \frac{1}{p_j^s} \frac{1}{p_k^s} \frac{1}{p_l^s} \frac{1}{p_m^s} \frac{1}{p_n^s}$$

\begin{center}
\bf.
\end{center}
\begin{center}
\bf.
\end{center}
\begin{center}
\bf.
\end{center}

\noindent The remainder function can also be defined as

\begin{equation}
{1+\displaystyle\sum\limits_{i=1}^\infty (-1)^{i}\displaystyle\sum\limits_{p}^{}{}_{i+1}}=\Big({\displaystyle\prod\limits_{k=1}^\infty {1-\frac{1}{\displaystyle p_k^s}•}}\Big)\Big({\displaystyle 1 + \displaystyle\sum\limits_{k=1}^\infty \frac{1}{p_k^s}\Big)} 
\end{equation}

\noindent which states that one plus the remainder function equals the product of the prime product function times the quantity one plus the prime zeta function \cite{PrimeZetaFunction}.\\

\noindent If we let the zeta function represent the prime product function found in the denominator of Eq. 2, we can make that substitution and rearrange Eq. 2 to obtain

\begin{equation}
{\displaystyle \zeta(s)}= \frac{\displaystyle 1 + \displaystyle\sum\limits_{k=1}^\infty \frac{1}{p_k^s}}{1+\displaystyle\sum\limits_{i=1}^\infty (-1)^{i}\displaystyle\sum\limits_{p}^{}{}_{i+1}}.
\end{equation}\\

\noindent In order to evaluate this new zeta function in the critical strip for complex $s$, we must employ the analytic version of the zeta function
 
\begin{equation}
\Big( 1-\frac{2}{2^s} \Big)^{-1}{\displaystyle\sum\limits_{n=1}^\infty} (-1)^{n+1} \frac{1}{n^s}
\end{equation}

\noindent formulated by Riemann so that we now have

\begin{equation}
\Big( 1-\frac{2}{2^s} \Big)^{-1}{\displaystyle\sum\limits_{n=1}^\infty} (-1)^{n+1} \frac{1}{n^s}= \frac{\displaystyle 1 + \displaystyle\sum\limits_{k=1}^\infty \frac{1}{p_k^s}}{1+\displaystyle\sum\limits_{i=1}^\infty (-1)^{i}\displaystyle\sum\limits_{p}^{}{}_{i+1}}.
\end{equation}

\noindent We will also need to introduce an analytic version of the prime zeta function into the numerator of the right-hand side of Eq. 7 so that it also will not blow up for values of $s$ for which $0<Re(s)<1$ in the critical strip of the complex plane.  To that end, we make an initial attempt to derive such a function that will parallel the analytic version of the Riemann zeta function on the left-hand side of Eq. 7 and be analyzable in the critical strip:\\

\begin{equation}
\Big( 1-\frac{2}{2^s} \Big)^{-1}{\displaystyle\sum\limits_{n=1}^\infty} (-1)^{n+1} \frac{1}{n^s}= \frac{\displaystyle 1 + \displaystyle\sum\limits_{k=1}^\infty (-1)^{k+1} \frac{1}{p_k^s}+ \displaystyle\sum\limits_{k=1}^\infty 2^k \displaystyle\sum\limits_{n=1}^\infty (-1)^{n+1} \frac{1}{p_{2^k{n}}^s}}{1+\displaystyle\sum\limits_{i=1}^\infty (-1)^{i}\displaystyle\sum\limits_{p}^{}{}_{i+1}}.
\end{equation}\\

\noindent The prime zeta function 

\begin{equation}
\displaystyle\sum\limits_{k=1}^\infty \frac{1}{p_k^s}
\end{equation}\\

\noindent found in the numerator of the right-hand side of Eq. 7 is now represented by

\begin{equation}
\displaystyle\sum\limits_{k=1}^\infty (-1)^{k+1} \frac{1}{p_k^s}+\displaystyle\sum\limits_{k=1}^\infty 2^k \displaystyle\sum\limits_{n=1}^\infty (-1)^{n+1} \frac{1}{p_{2^k{n}}^s}
\end{equation}\\

\noindent in the numerator of the right-hand side of Eq. 8.\\

\noindent However, neither functions (9) nor (10) enable analytic continuation in the critical strip.  Function (9) is the prime zeta function itself which does not converge for $0<Re(s)<1$; and the equivalent function (10) also does not yield a solution that is analytically continuable in the critical strip even though it employs the alternating prime zeta function.  This is because there is not an obvious way to eliminate every other prime number in the prime zeta function without having a remainder to account for which cannot be absorbed into the alternating prime zeta function itself to yield it analytically continuable in the critical strip.  This comes as no surprise since there is not an identifiable pattern in the sequence of prime numbers.  In light of this impasse, another prime number series is proposed to replace the prime zeta function in the numerator on the right-hand side of Eq. 7 which will enable an analytic continuation of the right-hand side of that equation in the critical strip:\\

\begin{equation}
\Big( 1-\frac{2}{2^s} \Big)^{-1}{\displaystyle\sum\limits_{n=1}^\infty} (-1)^{n+1} \frac{1}{n^s}= \frac{\displaystyle 1 + \Big( 1-\frac{2}{b^s} \Big)^{-1} \displaystyle\sum\limits_{k=1}^\infty (-1)^{k+1} \frac{1}{p_k^s}}{1+\displaystyle\sum\limits_{i=1}^\infty (-1)^{i}\displaystyle\sum\limits_{p}^{}{}_{i+1}}.
\end{equation}

\noindent On the right-hand side of Eq. 11, a variable $b$ is introduced to enable an analytic continuation of the alternating prime zeta function on both sides of the critical line $Re(s)=\frac{1}{2}$ within the critical strip.  We will refer to this complex variable as $beta$.\\

\noindent Referencing an earlier definition of the remainder function in Eq. 4, we now make a final substitution in Eq. 11 to arrive at a new definition of the zeta function which is totally analyzable in the critical strip and which is a function of prime sum and prime product functions exclusively:\\

\begin{equation}
\Big( 1-\frac{2}{2^s} \Big)^{-1}{\displaystyle\sum\limits_{n=1}^\infty} (-1)^{n+1} \frac{1}{n^s}= \frac{\displaystyle 1 + \Big( 1-\frac{2}{b^s} \Big)^{-1} \displaystyle\sum\limits_{k=1}^\infty (-1)^{k+1} \frac{1}{p_k^s}}{\Big({\displaystyle\prod\limits_{k=1}^\infty {1-\frac{1}{\displaystyle p_k^s}•}}\Big)\Big({\displaystyle 1 + \displaystyle\sum\limits_{k=1}^\infty \frac{1}{p_k^s}\Big)}}.
\end{equation}\\

\noindent We will refer to the right-hand side of Eq. 12 as the $meta$ function.  We can now solve Eq. 12 for the complex variable $beta$ of the alternating prime zeta function in terms of the analytic zeta function and the other prime number functions.  Moreover, we will also want to set the right-hand side of Eq. 12 equal to zero (since we are searching for the non-trivial zeros of the Riemann zeta function) and again solve for $beta$ so that we will have two equations for the solution of $b$ which we will be able to use simultaneously to find the zeros of the Riemann zeta function in the critical strip (i.e., we will have two equations in two unknowns).  It is worth noting at this point that since the solution for $beta$ will be a complex number, the evaluation of $b^s$ will be a complex number taken to a complex power.\\

\noindent The first solution for $beta$, obtained by equating the left-hand and right-hand sides of Eq. 12, is:\\

\begin{equation}
\mathlarger{\mathlarger{\mathlarger{b_{\zeta}=}}} \left\lbrace{\left\lbrace \frac{1}{2} \left\lbrace {1- {\frac {\displaystyle\sum\limits_{k=1}^\infty (-1)^{k+1} \frac{1}{p_k^s}} {\displaystyle \Big( 1-\frac{2}{2^s} \Big)^{-1} \displaystyle \sum \limits_{n=1}^\infty (-1)^{n+1} \frac{1}{n^s} \Big( {1 + \displaystyle\sum\limits_{k=1}^\infty \frac{1}{p_k^s} \Big) \Big( \displaystyle\prod\limits_{k=1}^\infty {1-\frac{1}{\displaystyle p_k^s}}} \Big) -1 }} }\right\rbrace\right\rbrace ^{-1}} \right\rbrace^{1/s}.
\end{equation}\\

\noindent The second solution of $beta$, obtained by setting the right-hand side of Eq. 12 equal to zero, is\\

\begin{equation}
\mathlarger{\mathlarger{\mathlarger{b_{m}=}}}
\left\lbrace \left\lbrace \frac{1}{2} \left\lbrace {1+\displaystyle\sum\limits_{k=1}^\infty (-1)^{k+1} \frac{1}{p_k^s}} \right\rbrace \right\rbrace ^{-1} \right\rbrace^{1/s}.
\end{equation}\\

\noindent For complex input $s$ where the $Re(s)>0$, the real parts of the two $b$ functions in Eqs. 13 and 14 are equal to each other, and the imaginary parts of the two $b$ functions are equal to each other, simultaneously, only at the locations of the non-trivial zeros of the Riemann zeta function located on the critical line $Re(s)=\frac{1}{2}$.  Conversely, when a complex input $s$ is not of the form $s=\frac{1}{2}+\rho i$, where $\rho$ is a non-trivial zero of the Riemann zeta function located on the critical line, then the real parts and the imaginary parts of the $beta$ functions in Eqs. 13 and 14, respectfully, will not equate simultaneously.  Thus, the zeros of the Riemann zeta function can be pinpointed on the critical line when the following two conditions are met:

\begin{equation}
\mathlarger{\mathlarger{\mathlarger{\text{Re}(b_{\zeta})}}}\mathlarger{\mathlarger{\mathlarger{=}}}\mathlarger{\mathlarger{\mathlarger{\text{Re}(b_{m})}}}
\end{equation}

\begin{center}
\underline{and}
\end{center}

\begin{equation}
\mathlarger{\mathlarger{\mathlarger{\text{Im}(b_{\zeta})}}}
\mathlarger{\mathlarger{\mathlarger{=}}}\mathlarger{\mathlarger{\mathlarger{\text{Im}(b_{m})}}}.
\end{equation}

\noindent The Mathematica \cite{Mathematica} plots in Appendix B provide strong empirical evidence that the conditions stated in Eqs. 15 and 16 are only met when the complex inputs $s$ are equal to the non-trivial zeros of the Riemann zeta function located on the critical line.  It is noted here that the exceptionally good degree of accuracy manifested in the Mathematica plots for $b_\zeta$ versus $b_m$ was obtained using a domain of only the first thousand prime numbers as inputs into the prime sum and prime product functions.  But while the plots in the appendices reveal some exciting properties of the $\zeta$, $m$ and $b$ functions, it is acknowledged that plots, in and of themselves, do not provide a validation of the Riemann Hypothesis.\\

\noindent However, a validation of RH is proposed in this paper by evidence of the convergence and divergence of the $b$, $\zeta$, and $m$ functions about the critical line $Re(s)=\frac{1}{2}$ within the critical strip of the complex plane.  When the imaginary part of $s$ is zero, the two solutions for $beta$ coincide exactly on the left-hand side of the critical strip but begin to diverge for $Re(s)>\frac{1}{2}$:

\begin{equation}
\mathlarger{\mathlarger{\mathlarger{b_{\zeta}=b_{m}}}} \,\,\,\,\,\,\,\,\,\,\,\, when \,\,\, Re(s)<\frac{1}{2} \,\,\, and \,\,\, Im(s)=0.
\end{equation}

\noindent Please see Figs. 1 and 2 for a visual depiction of this behavior \cite{Excel}.  For complex $s$ where the imaginary part of $s$ is not equal to zero, it was discovered that $Re(b_\zeta)=Re(b_m)$ and $Im(b_\zeta)=Im(b_m)$ only when $s=\frac{1}{2}+\rho i$, where $\rho$ is the imaginary part of any non-trivial zero of the Riemann zeta function.\\

\noindent Conversely, it was found that when the imaginary part of $s$ is zero, the $zeta$ and $meta$ functions coincide exactly on the right-hand side of the critical strip but begin to diverge for $Re(s)<\frac{1}{2}$:\\

\begin{equation}
\mathlarger{\mathlarger{\mathlarger{\zeta(s)=}}} \,\,\, \mathlarger{\mathlarger{\mathlarger{m(s)}}} \,\,\,\,\,\,\,\,\,\,\,\, when \,\,\, Re(s)>\frac{1}{2} \,\,\, and \,\,\, Im(s)=0.
\end{equation}\

\noindent Please refer to Figs. 1 and 2 for a visual depiction of this phenomenon \cite{Excel}.  For the case of complex $s$ where the imaginary part of $s$ is not equal to zero, $\zeta(s)=m(s)=0$ on the complex plane if and only if the complex variable $s$ is equal to $\frac{1}{2}+\rho i$, where $\rho$ is the imaginary part of any non-trivial zero of the Riemann zeta function.\\

\noindent Therefore, the only common location within the critical strip where all equalities are satisfied between the real and imaginary parts (respectfully) of the $beta$ functions, and where the $\zeta$ and $m$ functions are both equal to zero, are on the critical line $Re(s)=\frac{1}{2}$ at the locations of the non-trivial zeros of the Riemann zeta function (i.e., when $\zeta=meta=0$).  In contrast to the analytic version of the zeta function found on the left-hand side of Eq. 12, the new meta function on the right-hand side of Eq. 12 can directly be set equal to zero to locate the non-trivial zeros of the Riemann Hypothesis, and this will occur when the analytic version of the alternating prime zeta function found in the numerator on the right-hand side of Eq. 12 is equal to $-1$.  Thus, the new $meta$ function on the right-hand side of Eq. 12 enables one to find the non-trivial zeros of the Riemann zeta function directly by setting the right-hand side of that equation equal to zero.  To summarize, the condition \\\\

\begin{equation}
\mathlarger{\mathlarger{\mathlarger{\zeta(s)=}}} \,\,\, \mathlarger{\mathlarger{\mathlarger{m(s)}}} \,\,\, \mathlarger{\mathlarger{\mathlarger{=0}}}
\end{equation}\

\begin{center}
defined as
\end{center}

\begin{equation}
\Big( 1-\frac{2}{2^s} \Big)^{-1}{\displaystyle\sum\limits_{n=1}^\infty} (-1)^{n+1} \frac{1}{n^s}\,\,\,=\,\,\, \frac{\displaystyle 1 + \Big( 1-\frac{2}{b^s} \Big)^{-1} \displaystyle\sum\limits_{k=1}^\infty (-1)^{k+1} \frac{1}{p_k^s}}{\Big({\displaystyle\prod\limits_{k=1}^\infty {1-\frac{1}{\displaystyle p_k^s}•}}\Big)\Big({\displaystyle 1 + \displaystyle\sum\limits_{k=1}^\infty \frac{1}{p_k^s}\Big)}}\,\,\,=\,\,\,0
\end{equation}\\

\begin{center}
and the condition
\end{center}

\begin{equation}
\mathlarger{\mathlarger{\mathlarger{\text{Re}(b_{\zeta})=\text{Re}(b_{m})}}}\,\,\,\,\,\, and \,\,\,\,\,\,\mathlarger{\mathlarger{\mathlarger{\text{I}m(b_{\zeta})=\text{Im}(b_{m})}}}
\end{equation}\\

\noindent will only be met when the complex input $s$ is equal to a non-trivial zero of the Riemann function of the form $s=\frac{1}{2}+\rho i$.  Therefore, it is concluded that all non-trivial zeros of the Riemann zeta function must lie on the $Re(\frac{1}{2})$ line in the critical strip of the complex plane.\\

\begin{center}
\includegraphics[scale=0.7]{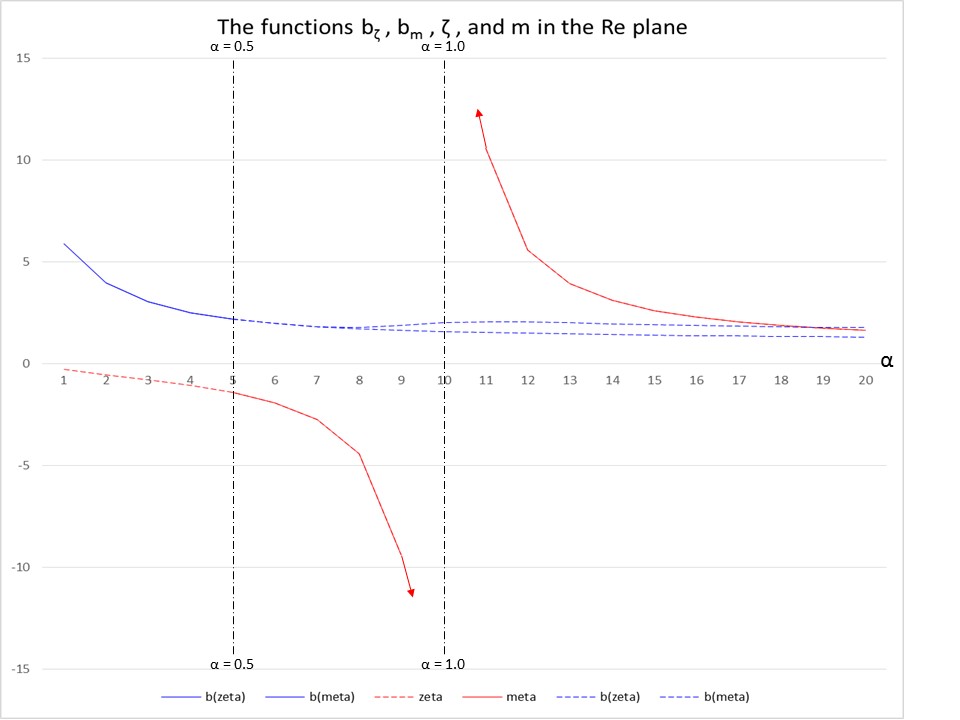}
\end{center}
\begin{center}
\textit{Fig. 1 - The beta functions and zeta/meta functions in the right-half of the real plane. Solid lines indicate where the functions converge, and the dotted lines indicate where the functions have diverged.}  
\end{center}

\begin{center}
\includegraphics[scale=0.7]{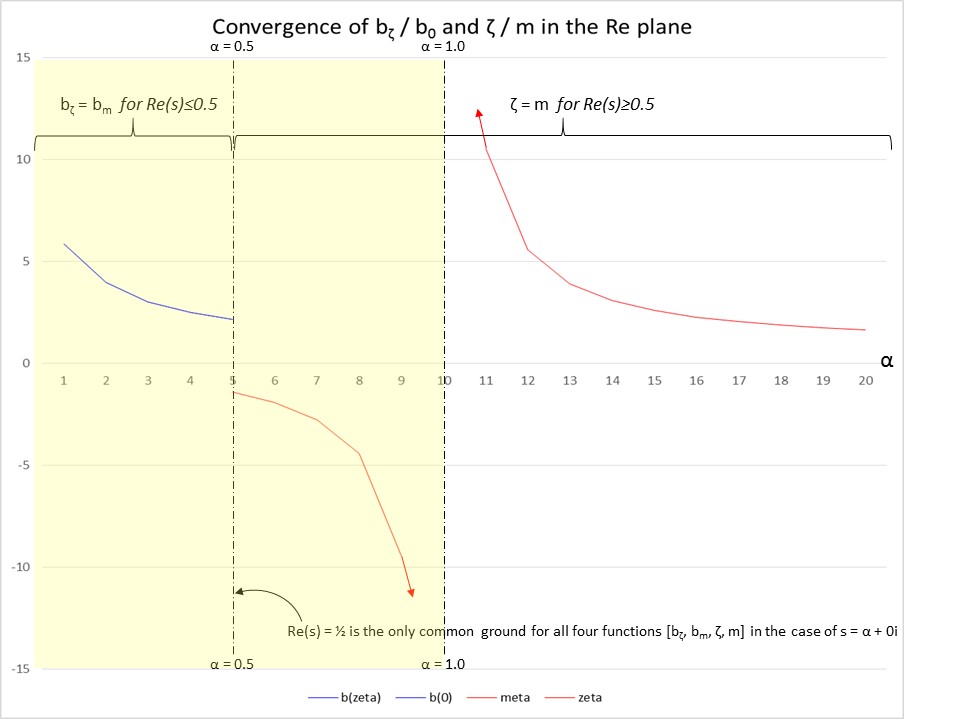}
\end{center}
\begin{center}
\textit{Fig. 2 - Convergence of the beta functions and the zeta/meta functions within the critical strip}  
\end{center}

\noindent To conclude, this discovery reveals a new method for locating zeros of the Riemann zeta function on the critical line using a computer algorithm and the formulae presented in this paper.  To do so, one would need to create a loop function which will pinpoint the imaginary zeros on the critical line to any desired degree of accuracy using the criteria $\text{Re}(b_{\zeta})=\text{Re}(b_{m})$ and $\text{I}m(b_{\zeta})=\text{Im}(b_{m})$.  Considerable time and complexity will be reduced in these computations of course by setting $Re(s)=\frac{1}{2}$.\\

\newpage
\appendix
\section{\\Plots of zeta vs. meta functions} \label{App:AppendixA}

This appendix displays combined plots of the $\zeta$ function versus the $m$ function for the first 50 non-trivial zeros of the Riemann zeta function.  Also included are plots of the $\zeta$ and $m$ functions with the following selected inputs of $s$:\\

$s=\alpha+0i$\\

$s=\alpha+14.134725...i$ (the imaginary part of the first zero of the Riemann zeta function)\\

$s=\alpha+143.111845...i$ (the imaginary part of the fiftieth zero of the Riemann zeta function)\\

\newpage
\begin{center}
\includegraphics[scale=0.67]{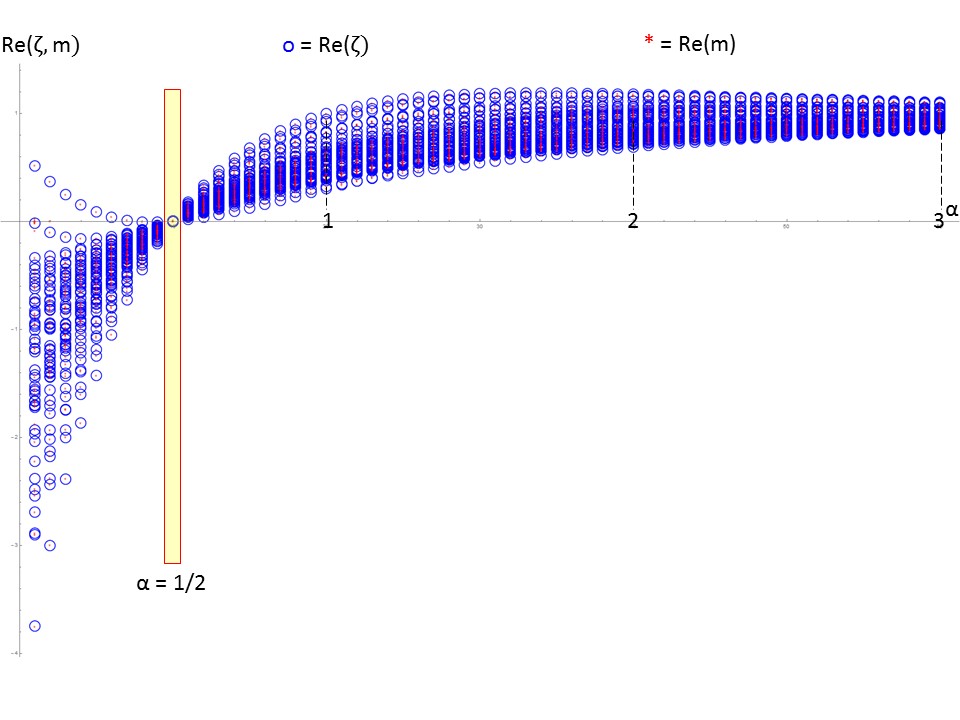}
\end{center}
\begin{center}
\textit{This is a plot of the real part of the zeta function $Re(\zeta)$ versus the real part of the meta function $Re(m)$ for the first 50 non-trivial zeros of the Riemann zeta function.  The yellow bar in this graph highlights the location of the convergence of the real parts of the zeta and meta functions, and this occurs only when both functions are equal to zero.  The real part of $\zeta$ and the real part of $m$ are both equal to zero only when the input $s$ is equal to $\frac{1}{2}+\rho i$, where $\rho$ is equal to the imaginary part of a non-trivial zero of the Riemann zeta function.  This result of the Riemann hypothesis is corroborated in this graph for the first 50 non-trivial zeros of the Riemann zeta function.}  
\end{center}

\newpage
\begin{center}
\includegraphics[scale=0.67]{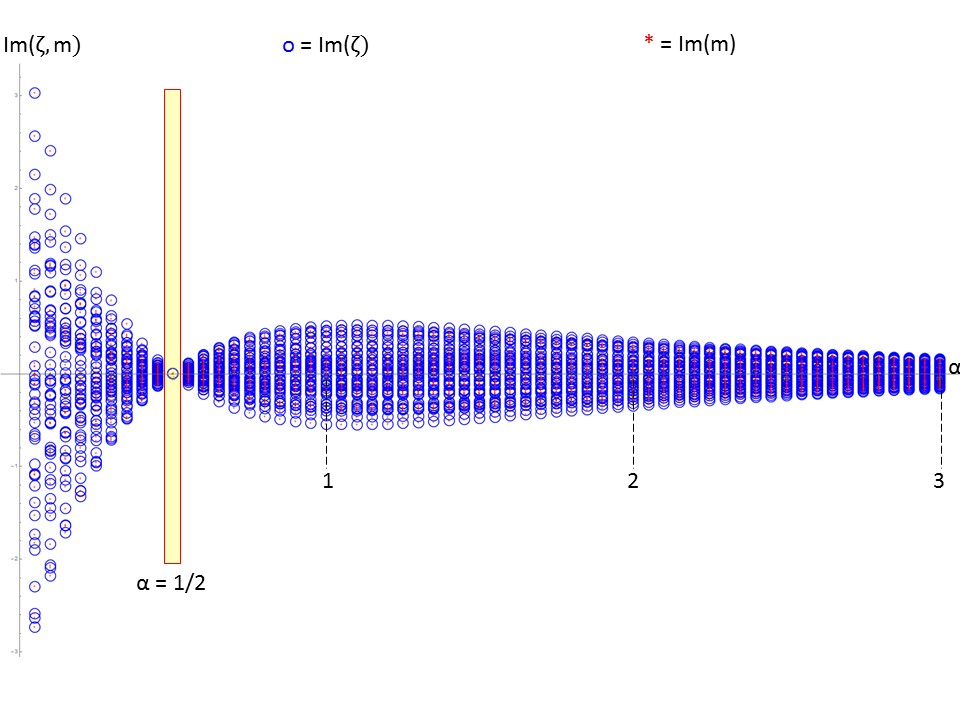}
\end{center}
\begin{center}
\textit{This is a plot of the imaginary part of the zeta function $Im(\zeta)$ versus the imaginary part of the meta function $Im(m)$ for the first 50 non-trivial zeros of the Riemann zeta function.  The yellow bar in this graph highlights the location of the convergence of the imaginary parts of the zeta and meta functions, and this occurs only when both functions are equal to zero.  Note that the imaginary part of $\zeta$ and the imaginary part of $m$ are both equal to zero only when the input $s$ is equal to $\frac{1}{2}+\rho i$, where $\rho$ is equal to the imaginary part of a non-trivial zero of the Riemann zeta function.  This result of the Riemann hypothesis is corroborated in this graph for the first 50 non-trivial zeros of the Riemann zeta function.}  
\end{center}

\newpage
\begin{center}
\includegraphics[scale=0.67]{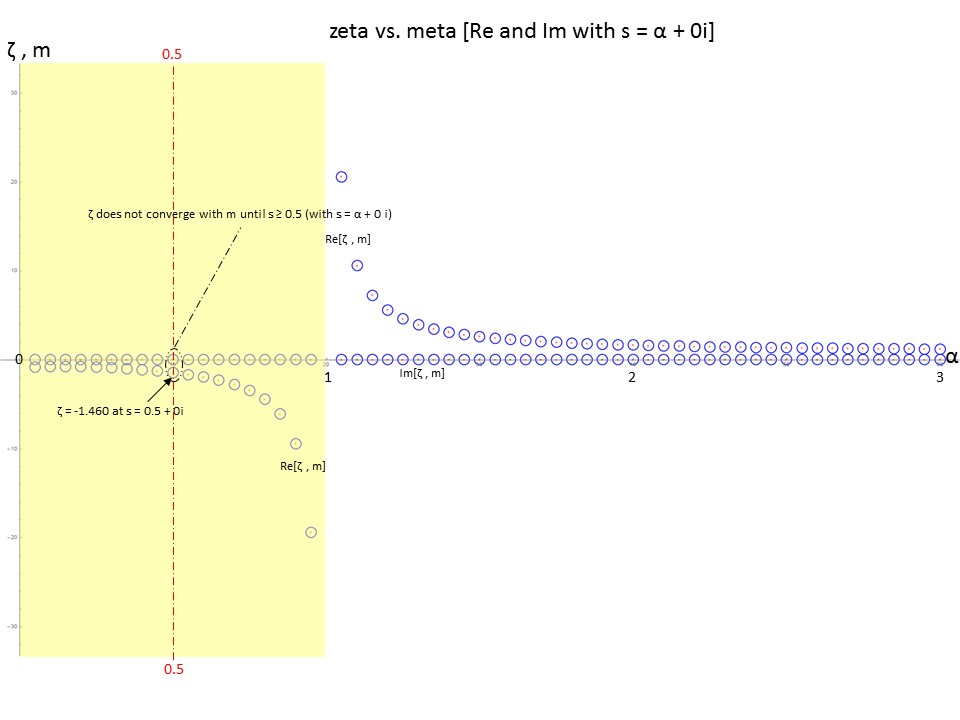}
\end{center}
\begin{center}
\textit{Plot of the real and imaginary parts of the zeta and meta functions for the real input $s=\alpha+0i$}  
\end{center}

\newpage
\begin{center}
\includegraphics[scale=0.67]{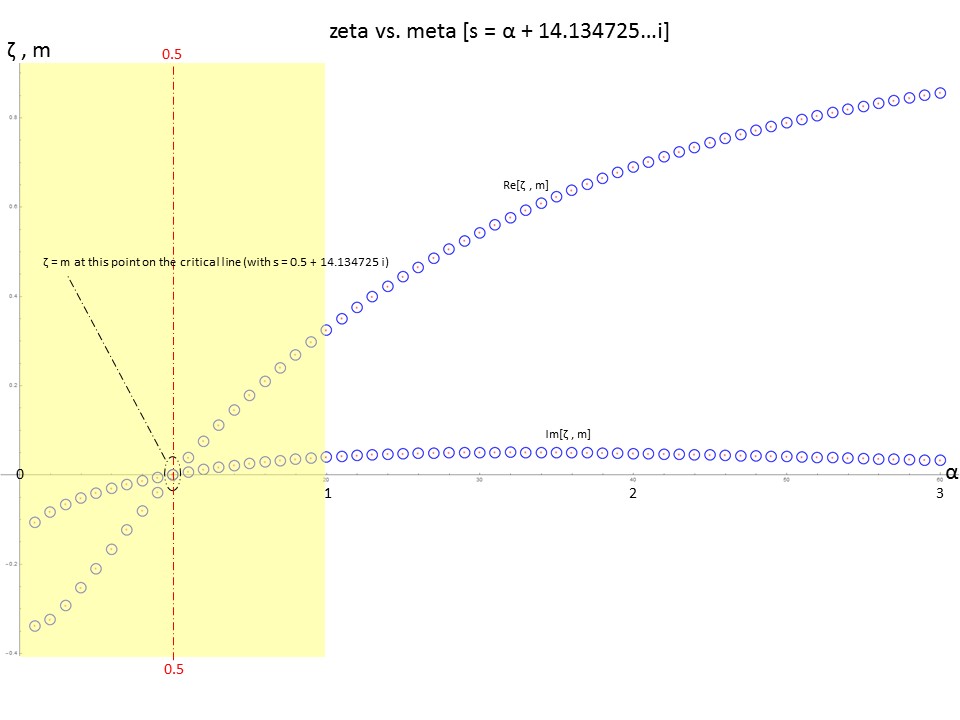}
\end{center}
\begin{center}
\textit{Plot of the real and imaginary parts of the zeta and meta functions for the complex input of $s=\alpha+14.134725i$}  
\end{center}

\newpage
\begin{center}
\includegraphics[scale=0.67]{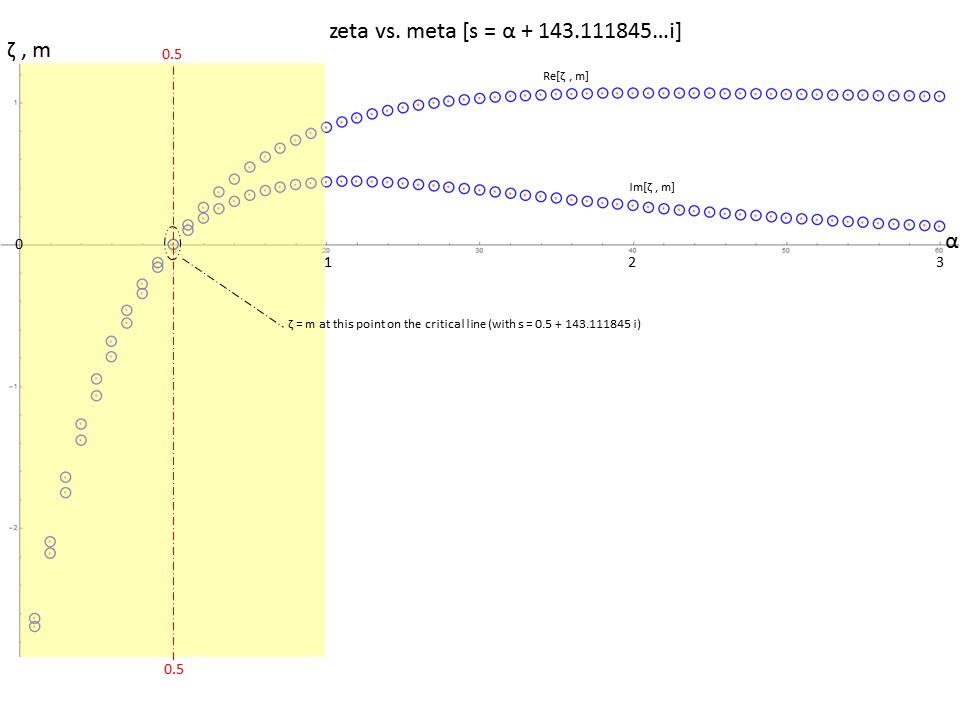}
\end{center}
\begin{center}
\textit{Plot of the real and imaginary parts of the zeta and meta functions for the complex input of $s=\alpha+143.111845i$}  
\end{center}

\newpage
\section{\\Plots of the $beta$ functions}
\label{App:AppendixB}

This appendix displays combined plots of the $b_\zeta$ function versus the $b_m$ function for the first 50 non-trivial zeros of the Riemann zeta function.  Also included are individual plots of the $b$ function for the real input $s=\alpha+0i$ and for the imaginary part of the first four non-trivial zeros of the Riemann zeta function, $s=\alpha+\rho i$.  The $b_\zeta$ function in all plots appears as blue circles and the $b_m$ function appears as red dots.

\newpage
\begin{center}
\includegraphics[scale=0.67]{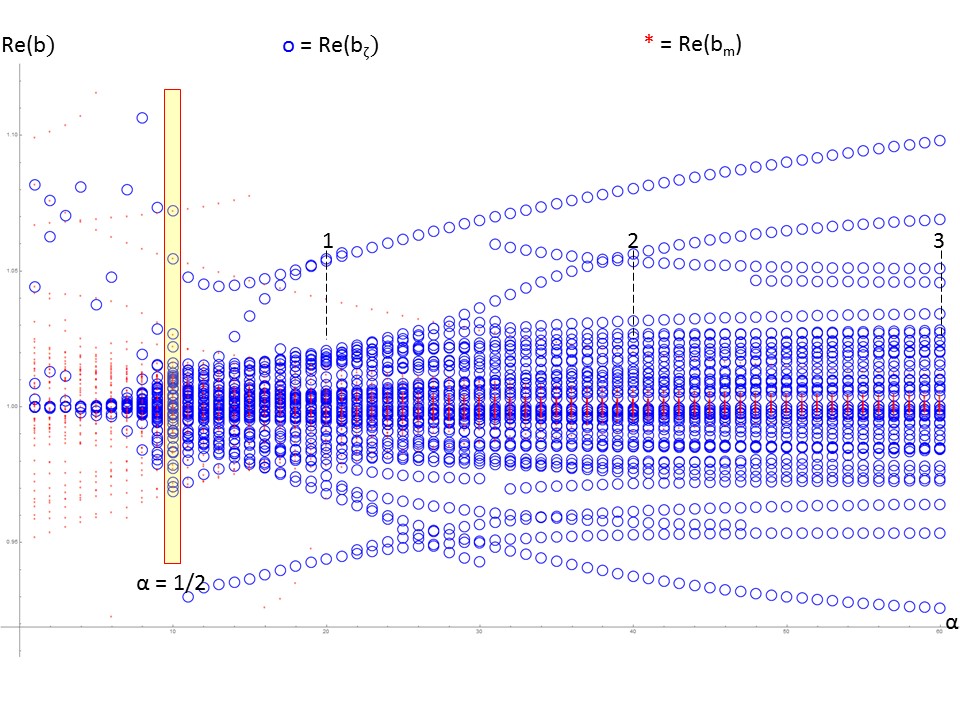}
\end{center}
\begin{center}
\textit{This is a plot of the real part of the $b_\zeta$ function $Re(b_\zeta)$ versus the real part of the $b_m$ function $Re(b_m)$ for the first 50 non-trivial zeros of the Riemann zeta function.  Note that the real part of $b_\zeta$ equals the real part of $b_m$ at every instance within the yellow band, i.e., where the input $s$ equals $\frac{1}{2}+\rho i$ where $\rho$ is the imaginary part of a non-trivial zero of the Riemann zeta function.  This is indicated by the red dots aligning to the center of the blue circles on the real $\frac{1}{2}$ line within the yellow band.  If instances of $Re(b_\zeta)$ and $Re(b_m)$ were added to this plot for any other input $\frac{1}{2}+it$ for which $t$ is not a non-trivial zero of the Riemann zeta function, then the two functions would not equate and the circles and dots would not align on the $\frac{1}{2}$ line within the yellow band.  It is interesting to note that in the midst of this seemingly chaotic graph, the real parts of $b_\zeta$ and $b_m$ coincide exactly for every instance of input $\frac{1}{2}+\rho i$ within the yellow band at $\alpha=\frac{1}{2}$.}  
\end{center}

\newpage
\begin{center}
\includegraphics[scale=0.67]{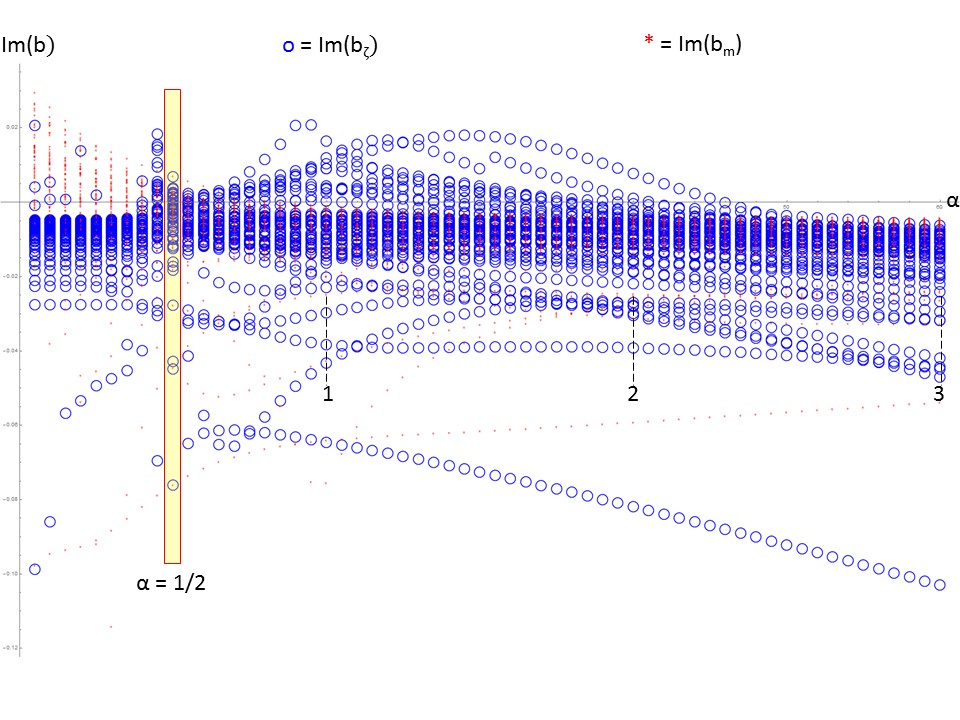}
\end{center}
\begin{center}
\textit{This is a plot of the imaginary part of the $b_\zeta$ function $Im(b_\zeta)$ versus the imaginary part of the $b_m$ function $Im(b_m)$ for the first 50 non-trivial zeros of the Riemann zeta function.  Note that the imaginary part of $b_\zeta$ equals the imaginary part of $b_m$ at every instance within the yellow band, i.e., where the input $s$ equals $\frac{1}{2}+\rho i$ where $\rho$ is the imaginary part of a non-trivial zero of the Riemann zeta function.  This is indicated by the red dots aligning to the center of the blue circles on the real $\frac{1}{2}$ line within the yellow band.  If instances of $Im(b_\zeta)$ and $Im(b_m)$ were added to this plot for any other input $\frac{1}{2}+it$ for which $t$ is not a non-trivial zero of the Riemann zeta function, then the two functions would not equate and the circles and dots would not align on the $\frac{1}{2}$ line within the yellow band.  It is interesting to note that in the midst of this seemingly chaotic graph, the imaginary parts of $b_\zeta$ and $b_m$ coincide exactly for every instance of input $\frac{1}{2}+\rho i$ within the yellow band at $\alpha=\frac{1}{2}$.}  
\end{center}

\newpage
\begin{center}
\includegraphics[scale=0.55]{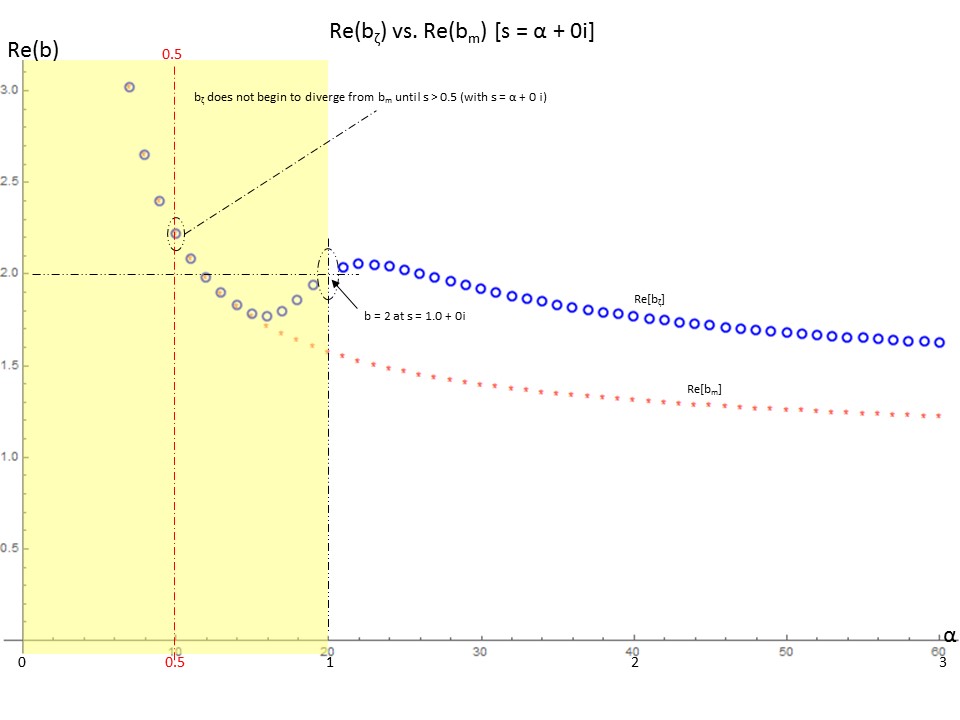}
\end{center}
\begin{center}
\textit{Plot of the real part of $b_\zeta$ versus the real part of $b_m$ for real input $s=\alpha+0i$}  
\end{center}

\begin{center}
\includegraphics[scale=0.55]{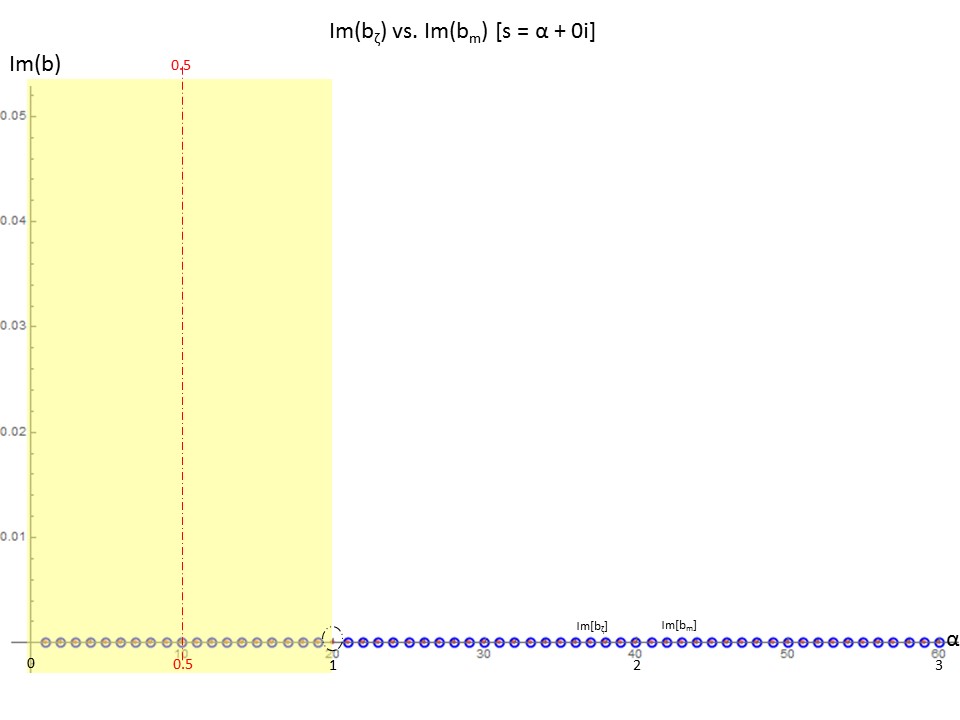}
\end{center}
\begin{center}
\textit{Plot of the imaginary part of $b_\zeta$ versus the imaginary part of $b_m$ for real input $s=\alpha+0i$}  
\end{center}

\newpage
\begin{center}
\includegraphics[scale=0.55]{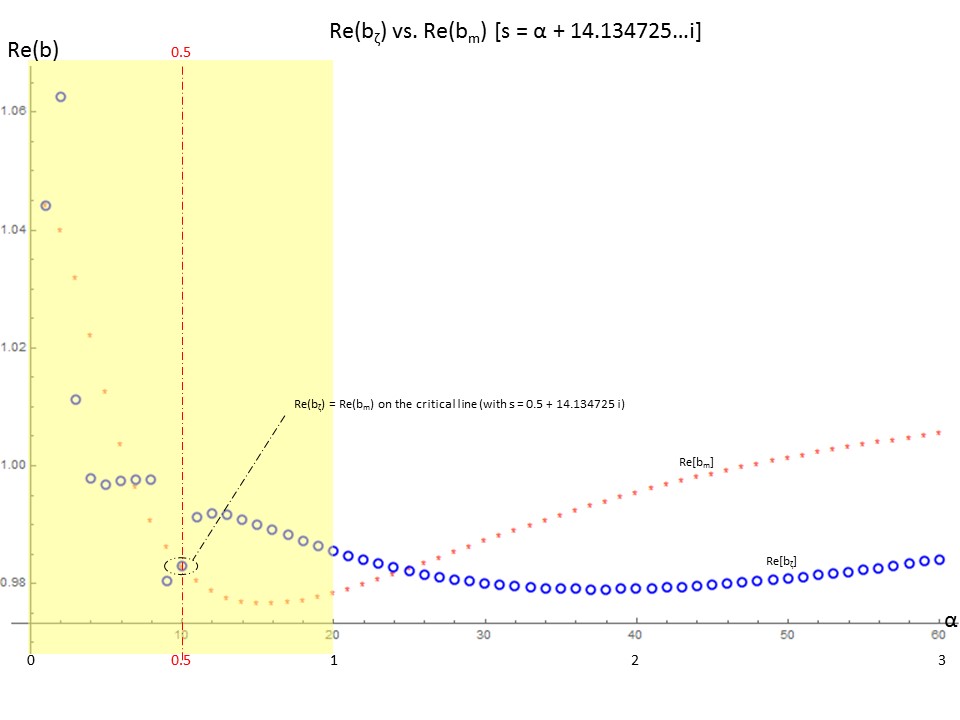}
\end{center}
\begin{center}
\textit{Plot of the real part of $b_\zeta$ versus the real part of $b_m$ for complex input $s=\alpha+14.134725...i$}  
\end{center}

\begin{center}
\includegraphics[scale=0.55]{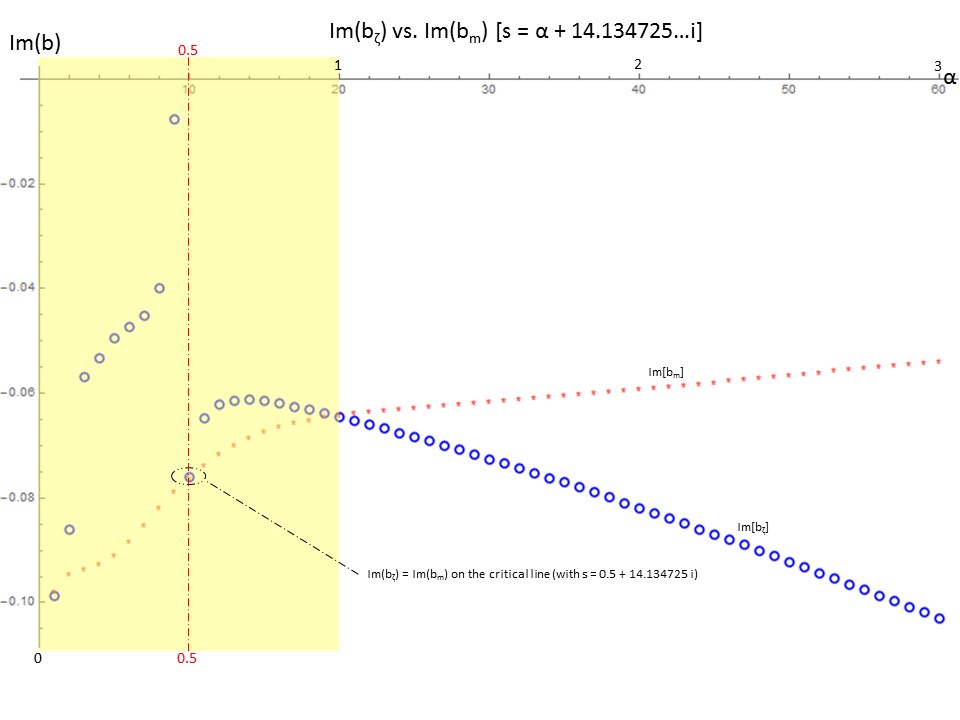}
\end{center}
\begin{center}
\textit{Plot of the imaginary part of $b_\zeta$ versus the imaginary part of $b_m$ for complex input $s=\alpha+14.134725...i$}  
\end{center}

\newpage
\begin{center}
\includegraphics[scale=0.55]{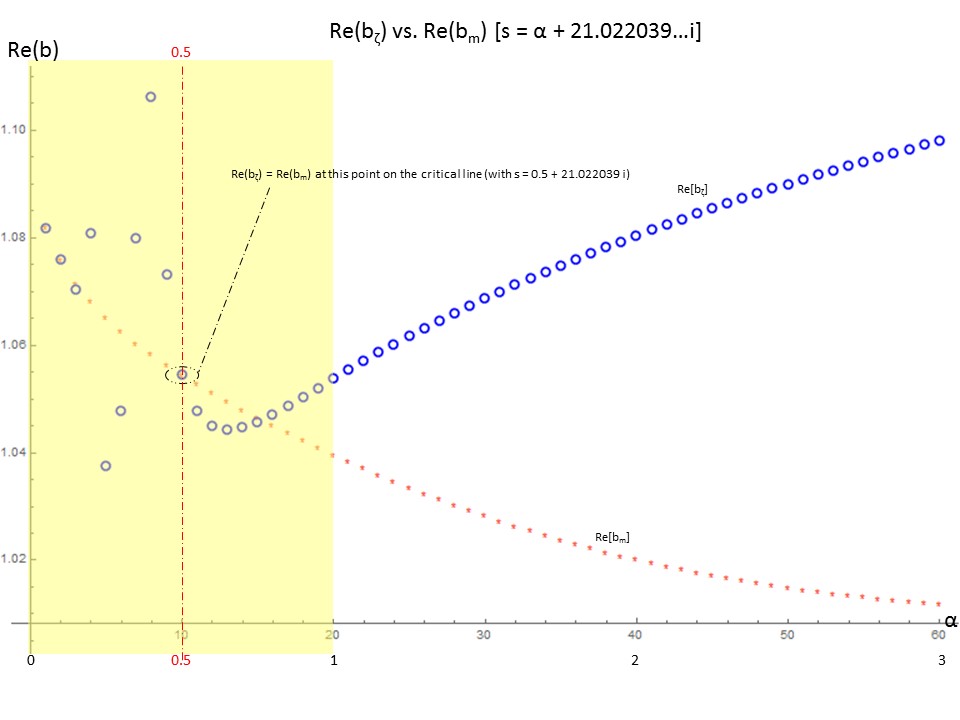}
\end{center}
\begin{center}
\textit{Plot of the real part of $b_\zeta$ versus the real part of $b_m$ for complex input $s=\alpha+21.022039...i$}  
\end{center}

\begin{center}
\includegraphics[scale=0.55]{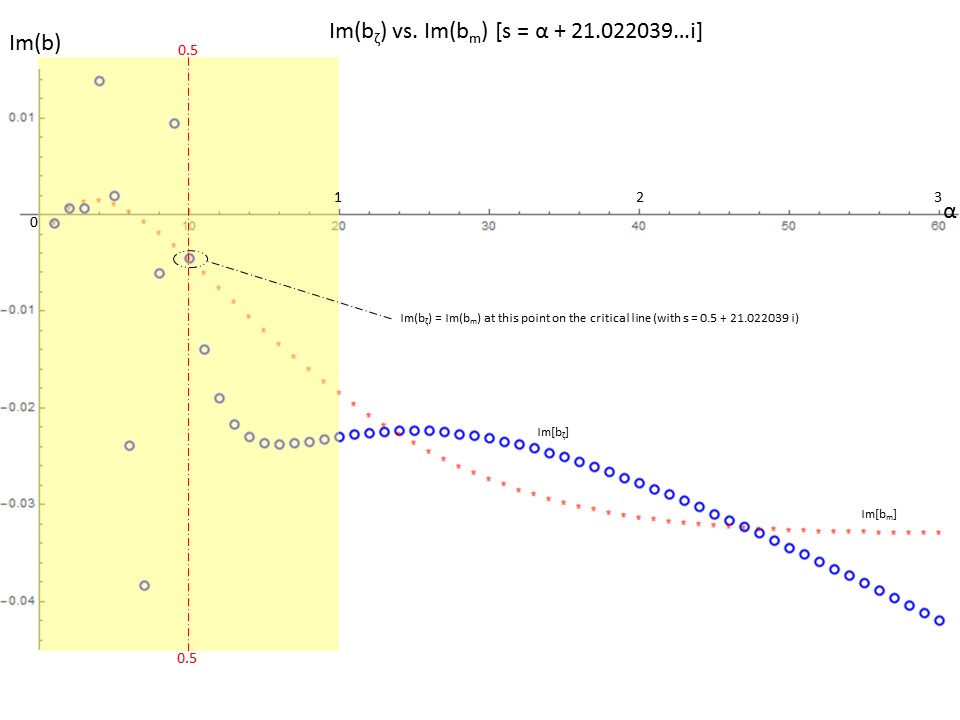}
\end{center}
\begin{center}
\textit{Plot of the imaginary part of $b_\zeta$ versus the imaginary part of $b_m$ for complex input $s=\alpha+21.022039...i$}  
\end{center}

\newpage
\begin{center}
\includegraphics[scale=0.55]{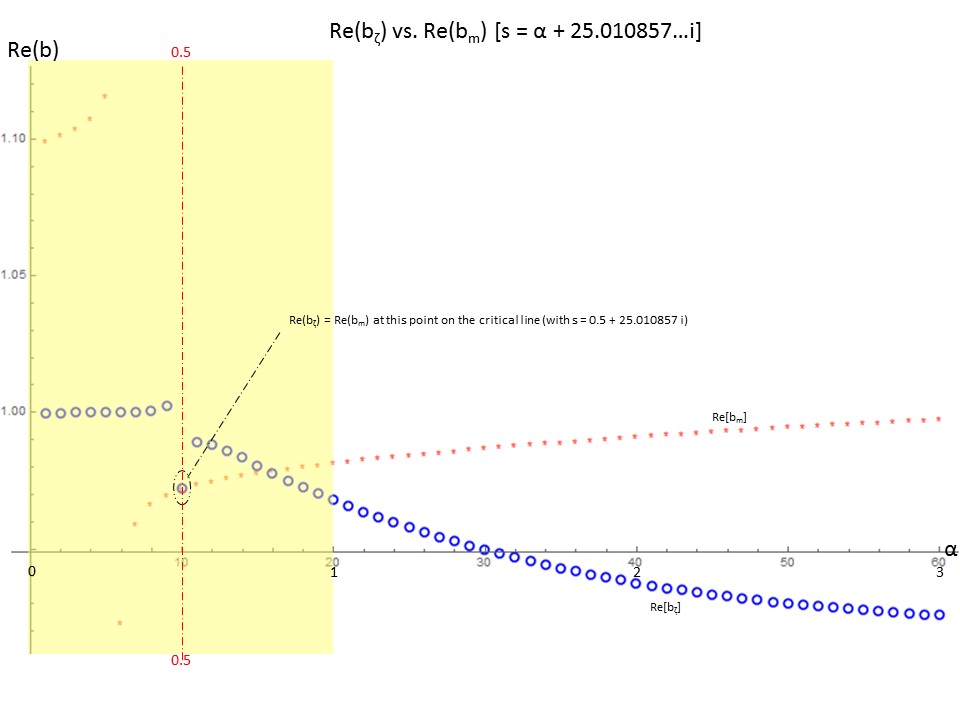}
\end{center}
\begin{center}
\textit{Plot of the real part of $b_\zeta$ versus the real part of $b_m$ for complex input $s=\alpha+25.010857...i$}  
\end{center}

\begin{center}
\includegraphics[scale=0.55]{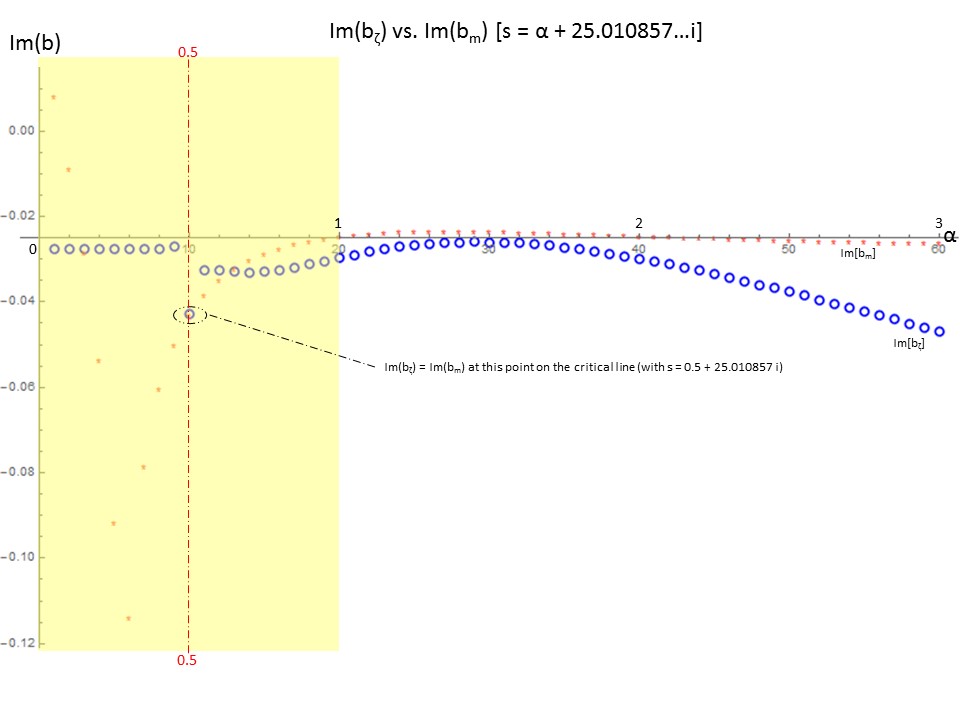}
\end{center}
\begin{center}
\textit{Plot of the imaginary part of $b_\zeta$ versus the imaginary part of $b_m$ for complex input $s=\alpha+25.010857...i$}  
\end{center}

\newpage
\begin{center}
\includegraphics[scale=0.55]{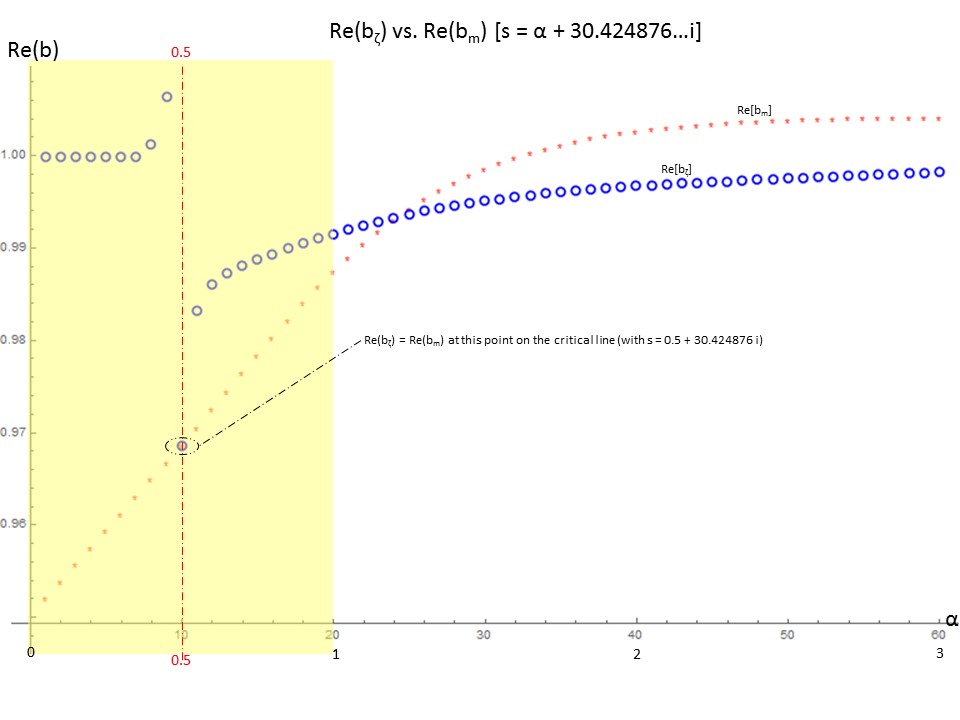}
\end{center}
\begin{center}
\textit{Plot of the real part of $b_\zeta$ versus the real part of $b_m$ for complex input $s=\alpha+30.424876...i$}  
\end{center}

\begin{center}
\includegraphics[scale=0.55]{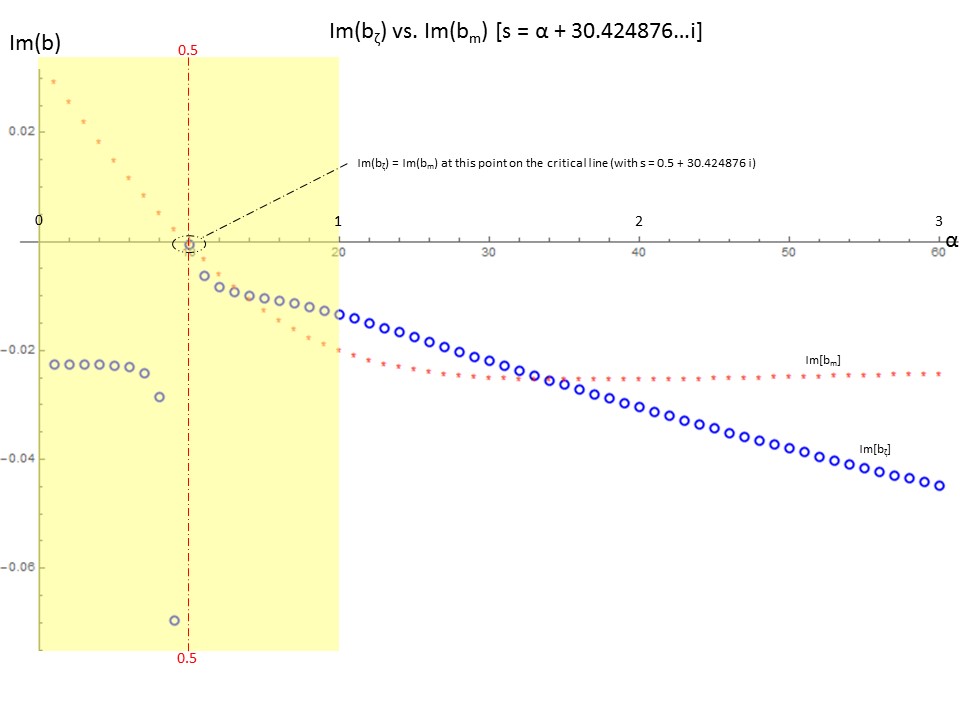}
\end{center}
\begin{center}
\textit{Plot of the imaginary part of $b_\zeta$ versus the imaginary part of $b_m$ for complex input $s=\alpha+30.424876...i$}  
\end{center}

\end{document}